\documentclass{amsart}
\usepackage{bm}

\def\AA{{\mathcal A}}

\def\CC{{\mathcal C}}
\def\DD{{\mathcal D}}
\def\EE{{\mathcal E}}
\def\FF{{\mathcal F}}
\def\GG{{\mathcal G}}
\def\HH{{\mathcal H}}

\def\KK{{\mathcal K}}
\def\LL{{\mathcal L}}
\def\MM{{\mathcal M}}

\def\RR{{\mathcal R}}

\def\XX{{\mathcal X}}

\def\bbR{\mathbb{R}}
\def\bbC{\mathbb{C}}
\def\bbZ{\mathbb{Z}}

\def\0{\mathbf{0}}

\def\qedrule{\hbox{\vrule height5.5pt depth0.5pt width3pt}}
\def\qed{\ifmmode \hskip 12pt plus 4pt minus 9pt \qedrule \else
\unskip \nobreak  \hskip 12pt plus 4pt minus 9pt \qedrule \medbreak \fi}

\def\<{\langle}
\def\>{\rangle}

\def\ga{\gamma}
\def\Ga{\Gamma}

\def\deg{\mathop{\rm deg}}

\def\Ref#1{(\ref{#1})}

\def\FT{\mathbb{F}}

\def\CB{\mathop{\rm CB}}
\def\HHH{\mathfrak{H}}

\def\MMM{\mathfrak{M}}
\def\TTT{\mathfrak{T}}

\def\phi{\varphi}

\newcommand\Ch{\gamma}

\newcommand\Domain{\bm{D}}
\newcommand\Range{\bm{R}}

\newtheorem{theorem}{Theorem}[section]
\newtheorem{proposition}[theorem]{Proposition}
\newtheorem{lemma}[theorem]{Lemma}
\newtheorem{corollary}[theorem]{Corollary}
\newtheorem*{theoremA}{Naimark's Theorem}

\theoremstyle{definition}
\newtheorem{remark}[theorem]{Remark}

\numberwithin{equation}{section}

\title{Extensions
of positive definite functions on free groups}

\author{M. Bakonyi and D. Timotin}
\address{Department of Mathematics and Statistics, Georgia State University, P.O. Box 4110, Atlanta, GA 30302-4110, USA}
\email{mbakonyi@gsu.edu}
\address{Institute of Mathematics of the Romanian Academy, PO Box 1-764, Bucharest 014700, Romania} \email{Dan.Timotin@imar.ro}
\date{}


\begin{document}


\begin{abstract}
An analogue of Krein's extension theorem is proved for operator-valued positive definite 
functions on free groups. The proof gives also the parametrization of all extensions by means of 
a generalized type of Szeg\"o parameters. One singles out a distinguished completion,
called central, which is related to quasi-multiplicative positive definite functions. An application is given
to factorization of noncommutative polynomials.
\end{abstract}

\maketitle

\section{Introduction}

Positive definite functions are an important object of
study in relation to group theory and $C^*$-algebras
since the basic work of Godement and Eymard~\cite{God, E}.
In the case of 
abelian groups, Bochner's theorem ensures that such a function is
the Fourier transform of a positive measure on the dual group, and
much of the theory develops along this line.
On nonabelian groups Fourier analysis is no more available, but
the focus is now on the relation to group representations. The
theory of positive definite functions is more intricate; for some notable
results, see also~\cite{Bo, BS, MF, Ha}.

The starting point of our investigation was the extension
problem for positive definite functions. The
earliest result seems to be Krein's extension theorem~\cite{K}, 
which says that every positive definite continuous
scalar function
on a real interval $(-a, a)$ admits a continuous
positive definite extension to ${\bbR}$; a recent generalization to totally
ordered groups can be found in~\cite{B1}. 

The analogue of Krein's theorem in $\bbZ^2$ is no more
true: a positive definite function defined on a 
rectangle symmetric with respect to the origin may
have no positive definite extension~\cite{Ru}. 
The sets with the property that every positive definite
function defined on them can be extended are characterized 
in~\cite{BN}.
A good reference for extension results (including
extensions from subgroups)
is~\cite[Chapter 4]{Sa}.

In the current paper we consider positive definite functions on the
most basic nonabelian group, namely the free group with an
arbitrary number of generators. (For general facts about Fourier
analysis on free groups, see~\cite{FT}.) 
We obtain an analogue of
Krein's theorem: a positive definite function
defined on the
set of words of length bounded by a fixed constant can be extended to the whole group. 
This is rather surprising in view of the result concerning $\bbZ^2$
discussed above; it is an illustration of a principle that has appeared
recently, namely that noncommutative objects might behave
better than their commutative analogues (see, for instance,~\cite{H}). 
Actually, our result is connected to~\cite{H} or~\cite{MC}, and in the
last section we deduce a factorization theorem for certain
polynomials in noncommutative variables. This is an area that has
received much interest in the last years \cite{MC, H, HM}.

The extension property is closely related to
a parametrization of all operator-valued positive definite functions
by means of sequences of contractions. These are an analogue of the \emph{choice sequences} of contractions, a generalization of the \emph{Szeg\"o parameters} in the
theory of orthogonal polynomials~\cite{Ge}, that has been developed by Foias and his coauthors in connection
with intertwining liftings (see~\cite{FF} 
and the references within). By using some
nontrivial graph theory, we are able to parallel the theory
of positive matrices as developed, for instance, in~\cite{C, C2}.

One should note that results in a close area of investigation
have been obtained by Popescu~\cite{Po0, Po1, Po2}. Most
notably, in~\cite{Po2} a similar extension problem is proved for the free
semigroup, together with a description of all solutions. The group case
that we have considered requires however
new arguments.

The plan of the paper is the following: in Section~\ref{se:prelim} we present preliminary
material concerning graphs and matrices. Section~\ref{se:grgr} introduces the free group,
and proves a basic result about its Cayley graph. The most important
part of the paper is Section~\ref{se:main}, where the main theorems concerning the structure
of positive definite functions and their extension properties are proved. In Section~\ref{se:central}
one discusses the existence of a distinguished extension, which is called the 
central one; it plays a role similar to that of central lifting in the theory
of intertwining liftings (or maximum entropy in data analysis)---see~\cite{BFF, FFG, FFGK}.
It is shown in Section~\ref{se:quasim} that some well known positive definite functions
on the free group may be obtined as central extensions. Finally,
in Section~\ref{se:noncom} we present the application to factorization of noncommutative
polynomials.

\section{Preliminaries}\label{se:prelim}

\subsection{Graphs}

A basic reference for graph theory is~\cite{GR}.
We consider undirected graphs $G=(V,E)$, where $V=V(G)$ denotes
the set of vertices and $E=E(G)$ the set of edges. If $\{v,w\}\in E(G)$, we
say that $v$ and $w$ are \emph{adjacent}.
$G$ is called \emph{complete} if any two vertices are adjacent.
 $G'=(V',E')$ is a \emph{subgraph} of $G$ if $V'\subset V$ and $E'\subset E$;
we write then $G'\subset G$.
If all edges of $G$ connecting two vertices in $V'$ are also edges of $G'$,
we say that $G'$ is the subgraph of $G$ \emph{induced} by $G'$; we
will write in this case $G'= G|_{V'}$. 
A set $C\subset V(G)$ is a \emph{clique} if the induced subgraph is complete.

A graph is 
called \emph{chordal} if
it contains no minimal cycles of length $>3$. 
An induced subgraph of a chordal graph is chordal.
We have the following basic result about chordal graphs~\cite{GJSW}.

\begin{lemma}\label{le:not2cliques}
If $G$ is chordal, then for any two nonadjacent vertices $v_1,v_2\in V(G)$ the set
of vertices adjacent both to $v_1$ and to $v_2$ is a clique.
\end{lemma}
%
%

A \emph{tree}  is a connected graph with
no cycles.
If $G$ is a tree, and $v,w\in V(G)$, there is a unique 
path $P(v,w)$ joining  $v$ and $w$ which
passes at most once through each vertex. We will call it the \emph{minimal path}
joining $v$ and $w$, and define
$d(v,w)$ to be its length (the number of edges it contains). 
For $n\ge 1$, we will denote by $\hat G_n$ the graph that has
the same vertices as $G$, while 
$E(\hat G_n)=\{(v,w) : d(v,w)\le n\}$ (in particular, $G=\hat G_1$).

\begin{lemma}\label{le:chordal}
If $G$ is a tree, then $\hat G_n$ is chordal for any $n\ge 1$.
\end{lemma}

\begin{proof} Take a minimal cycle $C$ of length $>3$ in $\hat G_n$. Suppose $x,y$
are elements of $C$ at a maximal distance. If $d(x,y)\le n$, then
$C$ is actually a clique, which is a contradiction. Thus $x$ and $y$ are
not adjacent in $\hat G_n$. Suppose $v,w$ are the two vertices of $\hat G_n$ adjacent
to $x$ in the cycle $C$. Now $P(x,v)$ has to pass through a vertex
which is on $P(x,y)$, since otherwise the union of these two paths
would be the minimal path connecting $y$ and $v$, and it would
have length strictly larger than $d(x,y)$. Denote by $v_0$ the
element of $P(x,v)\cap P(x,y)$ which has the largest distance to
$x$; since $d(y,v)=d(y,v_0)+d(v_0,v)\le d(y,x)=d(y,v_0)+d(v_0,x)$,
it follows that $d(v_0,v)\le d(v_0,x)$.

Similarly, if $w_0$ is the element of $P(x,w)\cap P(x,y)$ which has
the largest distance to $x$, it follows that $d(w_0,w)\le
d(w_0,x)$.

Suppose now that $d(v_0,x)\le d(w_0,x)$. Then
\begin{align*}
d(v,w)&= d(v,v_0)+d(v_0,w_0)+d(w_0,w)\\
&\le d(x, v_0)+d(v_0,w_0)+d(w_0,w)=d(x,w)\le n,
\end{align*}
since $w$ is adjacent to $x$. Then $(v,w)\in E$, and $C$ is not
minimal: a contradiction. Thus $\hat G_n$ is chordal.
\end{proof}

\begin{corollary}\label{co:chordal}
Suppose $G$ is a tree, and $x,y\in V(G)$, $d(x,y)=n+1$. Then the set 
\[
\{z\in V(G) :
\max(d(z,x), d(z,y))\le n\}
\] is a clique in $\hat G_n$.
\end{corollary}

\begin{proof}
We use Lemmas~\ref{le:not2cliques} and~\ref{le:chordal}.
\end{proof}

\subsection{Matrices and operators}

Suppose $A=(A_{ij})_{i,j\in I}$ is a block operator matrix, with 
entries $A_{ij}\in\LL(\HH)$, and $A\ge 0$ (as an operator on $\bigoplus_{i\in I} H$). There is an essentially unique 
Hilbert space, denoted by $\RR(A)$, together with isometries
$\omega_i(A): \HH\to \RR(A)$, such that $A_{ij}=\omega_i^*\omega_j$.
(If $\iota_i$ denotes the canonical embedding
of $\HH$ into the $i$th coordinate of $\bigoplus_{i\in I} \HH$, 
we may take, for instance, $\RR(A)$ to be the subspace 
of $\bigoplus_{i\in I} \HH$ spanned by all $\omega_i(A)\HH$, with
$\omega_i(A)=A^{1/2}\iota_i$.) 
We will call  $(\RR(A), \omega_i(A))$ a \emph{realization}
of $A$. 

If $J\subset I$, and $A'=A|J$ is the submatrix of $A$ obtained by taking only
rows and columns in~$J$, then one can embed isometrically
$\RR(A')$ into $\RR(A)$, 
such that $\omega_i(A)$ is just $\omega_i(A')$ followed by this embedding.


Suppose now that we are given only certain entries
$A_{ij}\in\LL(\HH)$ of a $p\times p$ block operator matrix, and
we are interested in \emph{completing} the remaining entries in
order to obtain a positive matrix. We assume that all $A_{ii}$ ($i\in I$)
are specified, and, since we are interested in
symmetric matrices, that
$A_{ij}$ and $A_{ji}$ are simultaneously specified.
The matrix is called \emph{partially positive} if all fully specified 
principal submatrices are
positive semidefinite. A similar definition applies to infinite
matrices.

The pattern of given
entries can be specified by a graph $G$, whose set of vertices is~$I$,
while $\{i,j\}$ ($i\not=j$) is an edge iff $A_{ij}$ is specified. 
An important role in completion
problems is then played by chordal graphs:
if $G$ is chordal, then every partially positive matrix with
pattern defined by $G$ has a positive semidefinite 
completion~\cite{GJSW}. 
A general reference for matrix completions is~\cite{J}.

In Section~4 we will use similar
arguments adapted to the group structure in order to extend positive definite
functions on the free group. But  we will only rely on the following  basic
result concerning completions by a single element (see, for instance,~\cite{ACC}
for a similar result).

\begin{lemma}\label{le:basic}
Suppose $A=(A_{ij})_{i,j\in I}$ is a partially positive block operator
matrix whose corresponding
pattern is the graph whose only missing edge is $\{k,l\}$. Then

(i) There exists a positive semidefinite completion of~$A$.

(ii) All such completions are in one-to-one correspondence with contractions
\[
\gamma:\RR(A|(I\setminus\{k\})\ominus \RR(A|(I\setminus\{k,l\})\to
\RR(A|(I\setminus\{l\})\ominus \RR(A|(I\setminus\{k,l\}).
\]
The correspondence is given by associating to each completion $\widetilde A$
the contraction 
\[
\gamma=P_{\RR(\widetilde A|(I\setminus\{l\})\ominus \RR(\widetilde A|(I\setminus\{k,l\})}|
\RR(\widetilde A|(I\setminus\{k\})\ominus \RR(\widetilde A|(I\setminus\{k,l\}).
\]
\end{lemma}

In particular, there is a distinguished completion, called \emph{central},
which corresponds to
$\gamma=0$.

\begin{remark}\label{re:details}
Although we intend to
avoid computational details related to matrix completion, 
below are sketched briefly some details concerning the correspondence 
stated in Lemma~\ref{le:basic}. Denote
\[
\begin{split}
\EE&=\RR(A|I\setminus\{k,l\}),\quad
\FF=\RR(A|I\setminus\{k\})\ominus \RR(A|I\setminus\{k,l\}),\\
\GG&=\RR(A|I\setminus\{l\})\ominus \RR(A|I\setminus\{k,l\}),
\end{split}
\]
and suppose $\omega_i:\HH\to \EE$, $i\in I\setminus\{k,l\}$
are the  embeddings defining the realization space of $A|I\setminus\{k,l\}$.
We may identify $\RR(A|I\setminus\{k\})$ with $\EE\oplus \FF$; with 
respect to this decomposition the embeddings $\omega'_i$
are given by
$\omega'_i= \begin{pmatrix} \omega_i\\ 0\end{pmatrix}$ for $i\not=l$, while
we will denote $\omega'_l= \begin{pmatrix} \alpha'\\ \beta'\end{pmatrix}$. 
Similarly, $\RR(A|I\setminus\{l\})$ is identified with $\EE\oplus \GG$, and
if $\omega''_i$ are the embeddings, then 
$\omega''_i= \begin{pmatrix} \omega_i\\ 0\end{pmatrix}$ for $i\not=k$, and
$\omega''_k= \begin{pmatrix} \alpha''\\ \beta''\end{pmatrix}$.

Suppose now that we have a contraction $\gamma:\FF\to\GG$, and denote
$D_\gamma=(I_\FF-\gamma^*\gamma)^{1/2}$, and $\DD_\gamma=\overline{D_\gamma\FF}$.
A concrete form for  the representing space of the associate
completion $\tilde A$ is then 
\begin{equation}\label{eq:defect1}
\RR(\tilde A)= \EE\oplus\GG \oplus \DD_{\gamma},
\end{equation}
where the embeddings $\tilde\omega_i:\HH\to\RR(\tilde A)$ are given by
\begin{equation}\label{eq:defect2}
\tilde\omega_i= \begin{pmatrix} \omega_i \\ 0 \\ 0  \end{pmatrix} \mbox{ for } i\not=
k, l, \quad
\tilde\omega_k= \begin{pmatrix}\alpha'' \\ \beta''\\ 0\end{pmatrix}, \quad
\tilde\omega_l= \begin{pmatrix}\alpha' \\ \gamma\beta'\\ D_\gamma\beta' \end{pmatrix}.
\end{equation}
Consequently, the formula for the completed entry is
\begin{equation}\label{eq:defect3}
A_{kl}=\tilde\omega^*_k \tilde\omega_l={\alpha''}^*\alpha' + {\beta''}^*\gamma\beta.
\end{equation}

Alternate formulas can be obtained by ``passing to the adjoint":
we have then $\RR(\tilde A)= \EE\oplus\FF \oplus \DD_{\gamma^*}$, with corresponding embeddings.
The essential uniqueness of the space $\RR(\tilde A)$ is reflected by a unitary 
$U_\gamma:\EE\oplus\GG \oplus \DD_{\gamma^*}\to  \EE\oplus\GG \oplus \DD_\gamma$
intertwining the two embeddings; this is given by
\begin{equation}\label{eq:ugamma}
U_\gamma = 
\begin{pmatrix}
I_\EE & 0 & 0\\
0 & \gamma & D_{\gamma^*}\\
0& D_\gamma & -\gamma^*
\end{pmatrix}.
\end{equation}
The operator appearing in the lower right corner is the well known ``Julia operator"~\cite{C}
related to $\gamma$.
\end{remark}


\section{Graphs and groups}\label{se:grgr}

We  consider in the sequel the group $\FT=\FT_m$, the free group
with $m$ generators $a_1,\dots,a_m$. Denote by $e$ the unit and by $\AA$
the set of generators of $\FT$.
Elements in $\FT$ are therefore words $s=b_1\dots b_n$, with letters
$b_i\in\AA\cup\AA^{-1}$; each word is usually written in its \emph{reduced} form, 
obtained after cancelling all products $aa^{-1}$. 
The length of a word $s$, denoted by $|s|$, is the
number of generators which appear in (the reduced form of) $s$. 
For a positive integer $n$,  define $S_n$ to be
the set of all words of length $\le n$ in $\FT$ and
$S'_n\subset S_n$ the set
of words of length exactly $n$; the number of 
elements in $S'_n$ is $2m (2m-1)^{n-1}$.

There is a notion that we will use repeatedly, and so we prefer
to introduce a notation. 
If $s,t$ are two words in reduced forms, we will call $u$ the  \emph{common beginning}
of $s$ and $t$ and write $u=\CB(s,t)$ if $s=us_1$ and $t=ut_1$, and $|t_1^{-1}s_1|=|s_1|+|t_1|$; in other
words, the first letters 
of $s_1$ and $t_1$ are different. For more than two words,
we define the common  beginning by induction:
\[
\CB(s_1,\dots s_p)=\CB(\CB(s_1,\dots,s_{p-1}), s_p).
\]
It can be seen that $\CB(s_1,\dots s_p)$ does not depend on the order of 
the elements $s_1,\dots s_p$, and is formed by their first common letters.

We include for completeness the proof of the following lemma.

\begin{lemma}\label{le:||square}
If $s\not=e$, then $|s^2|>|s|$; in particular, $s^2\not=e$.
\end{lemma}

\begin{proof}
When we write the word $s^2$, it may happen that some of the letters
at the beginning of $s$ cancel with some at the end of $s$. Let us
group the former in $s_1$ and the latter in $s_3$; thus $s=s_1 s_2 s_3$,
with $s_i$ in reduced form, $s_3s_1=e$, and the reduced form of $s^2$
is $s_1 s_2^2 s_3$. We have then $|s|=|s_1|+ |s_2|+ |s_3|$ and
$|s^2|=|s_1|+ 2|s_2|+ |s_3|$. If $|s^2|\le |s|$, we must have $s_2=e$. Therefore
$s=s_1s_3=e$.
\end{proof}

A total order $\preceq$ on $\FT$ will be called \emph{lexicographic} if,
for two elements $s_1,s_2\in\FT$, we have
 $s_1\preceq s_2$ whenever one of the following holds:

\begin{enumerate}
\item
$|s_1|<|s_2|$; 
\item
$|s_1|=|s_2|$, and, if $t=\CB(s_1,s_2)$, and $a_i$ is the first letter of $t^{-1}s_i$
(for $i=1,2$), then $a_1\preceq a_2$.
\end{enumerate}

We will write $s\prec t$ if $s\preceq t$ and $s\not=t$.
One can see that lexicographic orders are in one-to-one correspondence 
with their restrictions to 
$\AA\cup\AA^{-1}$. For the remainder of this section and for Section~\ref{se:main} a lexicographic
order $\preceq$ on $\FT$ will be fixed. 

Let $\sim$ be the equivalence relation on $\FT$ obtained by having the equivalence
classes $\hat s=\{s,s^{-1}\}$. We will denote also by $\preceq$ the order relation
on $\hat\FT=\FT/\sim$ defined by $\hat s\preceq \hat t$ iff $\min\{s,s^{-1}\}\preceq
\min\{t,t^{-1}\}$. If $\nu\in\hat\FT$, then $\nu^-$ and $\nu^+$ will be respectively
the predecessor and the successor of $\nu$ with respect to~$\preceq$.

\smallskip
There is a graph $\Ga$ naturally associated to $\FT$: the \emph{Cayley graph}
corresponding to the set of generators $\AA$.
Namely, $V(\Ga)$  are the elements of
$\FT$, while $s$ and $t$ are connected by an edge iff
$|s^{-1}t|=1$. Moreover, $\Ga$ is easily seen to be a tree: any
cycle would correspond to a nontrivial relation satisfied by the
generators of the group. The distance  $d$ between vertices of a tree defined in
Section~\ref{se:prelim} is $d(s,t)=|s^{-1}t|$. 

As a consequence of Lemma~\ref{le:chordal}, $\hat\Ga_n$ is chordal for any $n\ge1$. 
We will introduce a sequence of intermediate graphs $\Ga_\nu$, with $\nu\in\hat\FT$, as
follows. We have $V(\Ga_\nu)=\FT$ for all $\nu$, while 
$\{s,t\}\in E(\Ga_\nu)$ iff $\widehat{s^{-1}t}=\nu'$ for some $\nu'\preceq\nu$.
Obviously $\Ga_\nu\subset \Ga_{\nu'}$ for $\nu\preceq\nu'$, and $E(\Ga_{\nu^+})$
is obtained by adding to $E(\Ga_\nu)$ all edges $\{s,t\}$ with $\widehat{s^{-1}t}=\nu^+$.
Each $\Ga_\nu$ is invariant with respect to translations, and $\hat\Ga_n=\Ga_{\nu_n}$,
where $\nu_n$ is the last element in~$\hat S_n$.



The next proposition is the main technical ingredient of the paper.

\begin{proposition}\label{pr:2005}
With the above notation, $\Ga_\nu$ is chordal for all~$\nu$.
\end{proposition}

\begin{proof} If  $\Ga_\nu$ is not chordal, and $n=|\nu|$, we
may assume that $\nu$ is the last element in $\hat\FT$ of length $n$ with
this property.   Since $\hat\FT_n$ is chordal, $\nu\not=\nu_n$, 
and thus $|\nu^+|=n$.

Suppose then that 
$\Ga_\nu$ contains the cycle $(s_1,\dots,s_q)$, with $q\ge4$. At least
one of $\{s_1,s_3\}$ or $\{s_2,s_4\}$ must be an edge of $\Ga_{\nu^+}$, since
otherwise $\{s_1, s_2, s_3, s_4\}$ is a part of a cycle of length $\ge4$
in $\Ga_{\nu^+}$. We may assume that $\{s_1,s_3\} \in V(\Ga_{\nu^+})$; if we
denote $t=s_1^{-1}s_3$, then, since $\{s_1,s_3\} \not\in V(\Ga_\nu)$, we must
have $\hat t=\nu^+$.

Suppose that $q>4$. Then $\{s_1,s_4\}\not\in E(\Ga_\nu)$. If 
$\{s_1,s_4\} \in V(\Ga_{\nu^+})$, then $\widehat{s_1^{-1}s_4}=\nu^+$, and thus
either $s_1^{-1}s_4=t$ or $s_1^{-1}s_4=t^{-1}$. The first equality is 
impossible since it implies $s_3=s_4$. As for the second, it would lead
to $t^2=s_4^{-1}s_3$. But $|t|=|\nu^+|=n$, and thus, by Lemma~\ref{le:||square},
$|s_4^{-1}s_3|=|t^2|>n$. This contradicts $\{s_4,s_3\}\in V(\Ga_\nu)$;
consequently, we must have $q=4$.

Performing, if necessary, a translation, we may suppose
that the four cycle is $(e,s,t,r)$, and that $\{e,t\}\in E(\Ga_{\nu^+})$, 
(thus $t\in\nu^+$, and $t\prec t^{-1}$).
Thus, $e,s,r,t$ are all different,
$\{e,s\}, \{s,t\}, \{r,t\}, \{e,r\}$
are edges of $G_\nu$ while $\{e,t\}, \{s,r\}$ are not, 
and $\{e,t\}$ is an edge of $G_{\nu^+}$.
These assumptions imply that:
\begin{equation}\label{eq:ineq}
|t|=n,\quad |s|\le n,\quad |r|\le n, \quad |t^{-1}s|\le n, \quad |t^{-1}r|\le n,
\quad |r^{-1}s|\ge n.
\end{equation}

Let us denote $u=\CB(s,r,t)$. We
have $r=ur_1$, $t=ut_1$, $s=us_1$, and 
at least two among the elements
$\CB(s_1,r_1)$, $\CB(r_1,t_1)$, $\CB(s_1,t_1)$ are equal to $e$.

Suppose first that this happens with $(s_1,t_1)$ and $(r_1,t_1)$. If
$v=\CB(s_1,r_1)$, then $s_1=vs_2$, $r_1=vr_2$. The inequalities~\eqref{eq:ineq}
imply
\begin{align}
|u|+|t_1|&=n,\label{eq:32}\\
|u|+|v|+|s_2|&\le n, \quad |u|+|v|+|r_2|\le n,\label{eq:33}\\
|t_1|+|v|+|s_2|&\le n, \quad |t_1|+|v|+|r_2|\le n,\label{eq:35}\\
|s_2|+|r_2|&\ge n.\label{eq:34}
\end{align}
Replacing $|t_1|=n-|u|$ from~\eqref{eq:32}, one obtains from~\eqref{eq:35}
$|v|+|s_2|\le |u|$, $|v|+|r_2|\le |u|$. Then~\eqref{eq:33} 
imply $|v|+|s_2|\le n/2$, $|v|+|r_2|\le n/2$.
Comparing the last two inequalities with~\eqref{eq:34}, one obtains
$|v|=0$. This means that all pairs $(s_1,r_1), (r_1,t_1)$ and $(s_1,t_1)$ have as common beginning $e$, and so we may as well assume from the start that 
$(s_1,r_1)$ and $(s_1,t_1)$ have as common beginning $e$ (the case
 $(s_1,r_1)$ and  $(r_1,t_1)$ is symmetrical, and can be treated likewise; note
 that the whole situation is symmetric in $s$ and~$r$).
\begin{figure}[h]

\begin{picture}(600,150)(-60,10)
 \put(10,75){$\circ$}
 \put(10,85){$e$}
 \put(15,78){\line(1,0){20}}
 \put(38,76){\dots}
 \put(52,78){\line(1,0){20}}
 \put(75,75){$\circ$}
 \put(75,65){$u$}
 
 \put(77,83){\line(0,1){20}}
 \put(76,106){$\vdots$}
 \put(77,119){\line(0,1){20}}
 \put(75,140){$\circ$}
 \put(82,140){$s=us_1$}
 \put(85,110){$s_1$}
 
 \put(80,78){\line(1,0){20}}
 \put(103,76){\dots}
 \put(117,78){\line(1,0){20}}
 \put(140,75){$\circ$}
 \put(138,85){$uw$}
 
 \put(105,85){$w$}
 
 \put(142,18){\line(0,1){20}}
 \put(142,41){$\vdots$}
 \put(142,54){\line(0,1){20}}
 \put(140,10){$\circ$}
 \put(147,10){$r=uwr_2$}
 \put(151,41){$r_2$}
 
  \put(146,78){\line(1,0){20}}
 \put(169,76){\dots}
 \put(183,78){\line(1,0){20}}
 \put(206,75){$\circ$}
 \put(216,75){$t=uwt_2$}
 \put(172,85){$t_2$}
 
\end{picture}
\caption{}
\end{figure}

Then we will assume in the sequel that $w=\CB(t_1,r_1)$; thus $r_1=wr_2$ and 
$t_1=wt_2$ (see Figure~1). 
We have then 
\begin{equation}\label{eq:2}
|u|+|w|+|t_2|=n,\quad |u|+|w|+|r_2|\le n, \quad |s_1|+|w|+|t_2|\le n,
\end{equation}
whence $|s_1|\le |u|$, $|r_2|\le |t_2|$. But, since $\{s,r\}$ is not an
edge of $\Ga_\nu$, we have $|s_1|+|w|+|r_2|\ge n$. Comparing this last 
inequality with the first inequality in~\eqref{eq:2}, it follows
that
\[
|s_1|=|u|,\quad |r_2|=|t_2|.
\]

Now, $s^{-1}t$ is a word of length $n$ different from $t$;
thus $s_1^{-1}w\not=uw$, and
$|s_1^{-1}w|=|uw|$.

On the other hand, $t\in\nu^+$ begins with $uw$, and
$t\preceq t^{-1}$. Suppose  
$t'\in\FT$ with $|t'|=n$, begins with $t'_1$, with $|t'_1|=|uw|$. If  $t'_1\prec uw$, then 
the definition of the lexicographic order implies that
$t'\prec t$ and $\hat t'\prec \hat t$. Applying this argument
to $t'=s^{-1}r$, $t'_1=s_1^{-1}w$, it follows that, if $s_1^{-1}w\prec uw$,
then $\widehat{s^{-1}r}\prec \nu^+$. Then $\widehat{s^{-1}r}\prec \nu$,
and therefore $\{s,r\}\in E(\Ga_\nu)$: a contradiction.
Therefore $uw\prec s_1^{-1}w$, and $t\prec s^{-1}t$.

Since, however, $\{s,t\}\in\Ga_\nu$, it follows that $\widehat{t^{-1}s}\preceq\nu$, and thus 
\begin{equation}\label{eq:4}
t^{-1}s\preceq t\preceq t^{-1}.
\end{equation}
Also, $\{1,s\}\in E(\Ga_\nu)$ implies $|s|=|u|+|s_1|=2|u|\le n$, and thus 
$|u|\le n/2$. Therefore $|t_2^{-1}w^{-1}|=|wt_2|\ge n/2$.

Now $t_2^{-1}w^{-1}=\CB(t^{-1}s, t^{-1})$, and~\eqref{eq:4} implies
that $t_2^{-1}w^{-1}$ is also the beginning of $t$. Thus $\CB(t,t^{-1})$
has length $\ell\ge n/2$. Writing then $t^2=(t^{-1})^{-1}t$, $\CB(t,t^{-1})$ cancels, and we obtain $|t^2|\le 2n-2\ell\le n=|t|$. By
Lemma~\ref{le:||square}, this would imply $t=e$: a contradiction, since the elements $e,s,r,t$
of the assumed 4-cycle must be distinct. The proposition is proved.
\end{proof}



If $\nu\in\hat\FT$, then
there exist cliques in $\Ga_{\nu^+}$ which are not cliques in $\Ga_\nu$:
we may start with any edge of $\Ga_{\nu^+}$ which is not an edge of
$\Ga_\nu$ and take a maximal clique (in $\Ga_{\nu^+}$) which contains it.
Such a clique is necessarily finite, since the length of edges
is bounded by $|\nu^+|$.

\begin{corollary}\label{co:cliques}
Suppose $\nu\in\hat\FT$, and 
$C$ is a clique in $\Ga_{\nu^+}$ which is not a clique in $\Ga_\nu$.
Then $C$
contains a single edge in $\Ga_{\nu^+}$ which is not an edge in $\Ga_\nu$.
\end{corollary}

\begin{proof}
Suppose $C$ contains $\{s,t\}, \{s',t'\}$ with $s^{-1}t=s'{}^{-1}t'=v\in\nu^+$.
Obviously $s\not=s'$ and $t\not=t'$.
If $s=t'$, then $s'{}^{-1}t=v^2$, and
thus  the edge $\{s',t\}$ has length strictly larger than $n=|v|$ by Lemma~\ref{le:||square}. 
One shows similarly that $s'\not=t$. Then $(s,t,t',s')$ is a 4-cycle in $\Ga_\nu$, 
which contradicts 
Proposition~\ref{pr:2005}.
\end{proof}

\section{Positive definite functions}\label{se:main}

A function $\Phi:\FT\to\LL(\HH)$ is \emph{positive definite} if for every
$s_1,\dots, s_k\in\FT$  the operator matrix 
$A(\Phi;\{s_1,\dots, s_k\}):=[\Phi(s_i^{-1}s_j)]_{i,j=1}^k$ is positive semidefinite.

In general, for a finite set $S\subset\FT$, we will use the notation
$A(\Phi; S):= [\Phi(s^{-1}t)]_{s,t\in S}$. There is here a slight abuse
of notation, since the matrix in the right hand side of the equality 
depends on the order of the elements of $S$, and changing the order 
amounts to intertwining the rows and columns of $A(\Phi; S) $; however,
the reader can easily check that this ambiguity is irrelevant in all instances where 
this notation is used below.

Positive definite functions are connected to representations
of $\FT$. The relation is given by Naimark's Dilation Theorem~\cite{P, Sa}: 

\begin{theoremA}\label{th:naimark}
The function $\Phi:\FT\to\LL(\HH)$ is positive definite if and only if
there exists a representation $\pi$ of\/ $\FT$ on a Hilbert space
$\KK$, and an operator $V:\HH\to\KK$, such that, for all
$s,t\in\FT$, $\Phi(s^{-1}t)=V^*\pi(s^{-1}t)V$.
\end{theoremA}

Such a representation $\pi$ is called a \emph{dilation} of $\Phi$; it is uniquely
determined, up to unitary isomorphism, by the minimality condition
$\KK=\bigvee_{s\in\FT}\pi(s)V\HH$.

It is immediate from the definition of positive definiteness that
$\Phi(e)$ is a positive operator. A standard argument shows
that, if $\HH_0=\overline{\Phi(e)\HH}$, then there exists
a positive definite function $\Phi_0:\FT\to \LL(\HH_0)$, such
that $\Phi_0(e)=I_{\HH_0}$, and
$\Phi(s)h=\Phi_0(s)\Phi(e)^{1/2}h$ for all $h\in\HH$.
We will suppose in the sequel that $\Phi(e)=I_\HH$;
equivalently, the operator~$V$ in Naimark's Theorem
is an isometry.

We will also consider positive definite functions defined on subsets of $\FT$.
If $\Sigma\subset\FT$ such that $\Sigma=\Sigma^{-1}$, then a function
$\phi:\Sigma\to\LL(\HH)$ is called \emph{positive definite} if for every
$s_1,\dots, s_k\in\FT$ such that $s_i^{-1} s_j\in \Sigma$ for every
$i,j=1,\dots, k$, the operator matrix $A(\Phi;\{s_1,\dots, s_k\})$ is positive semidefinite.
Obviously, if $\Phi: \FT\to\LL(\HH)$ is positive definite, then
$\Phi|\Sigma$ is positive definite for all $\Sigma\subset\FT$.

Remember now that we have a fixed lexicographic order on $\FT$, and
consequently on $\hat\FT$.
Define  $\Sigma_\nu=\bigcup_{\nu'\preceq\nu}\nu$;  the following lemma is then
obvious, on closer inspection.

\begin{lemma}\label{le:obvious} Suppose $\mu\in\FT$, and $\phi:\Sigma_\mu\to\LL(\HH)$ 
is positive definite. Then 
the map $\Phi\mapsto (\Phi|\Sigma_\nu)_{\nu\in\hat\FT, \mu\preceq\nu}$ defines
a one-to-one correspondence between:

1) positive definite functions on $\FT$ whose
restriction to $\Sigma_\mu$ is $\phi$;

\noindent and

2) sequences $(\phi_\nu)_{\nu\in\hat\FT, \mu\preceq\nu}$, $\phi_\nu$ positive definite 
function defined on $\Sigma_\nu$, such that, if $\nu\preceq\nu'$,
then $\phi_\nu=\phi_{\nu'}|\Sigma_\nu$.
\end{lemma}

Note that the lemma may be applied to the case
$\mu=\{e\}$ and $\phi(e)=I_\HH$, when we obtain the
correspondence between a positive definite function on the whole
group $\FT$ and the sequence of its restrictions.

Let us now fix $\nu\in\hat\FT$, $\nu\not=\{e\}$; thus $\nu=\{s_\nu, s_\nu^{-1}\}$, 
with $s_\nu\preceq s_\nu^{-1}$. Suppose that $\phi_{\nu^-}$ is a positive
definite function on $\Sigma_{\nu^-}$.
Consider the maximal clique $C_\nu$ in $\Ga_\nu$ that contains $\{e, s_\nu\}$;
by Corollary~\ref{co:cliques} $\{e, s_\nu\}$  is the unique edge of $C_\nu$ that
is not in $E(\Ga_{\nu^-})$. Consequently, the following definitions make sense:
\begin{equation}\label{eq:domrange}
\begin{split}
\Domain(\phi_{\nu^-})&= \RR(A(\phi_{\nu^-};C_\nu\setminus\{s_\nu\}))
\ominus \RR(A(\phi_{\nu^-};C_\nu\setminus\{e,s_\nu\}))\\
\Range(\phi_{\nu^-})&= \RR(A(\phi_{\nu^-};C_\nu\setminus\{e\}))
\ominus \RR(A(\phi_{\nu^-};C_\nu\setminus\{e,s_\nu\}))
\end{split}
\end{equation}


\begin{lemma}\label{le:step}With the above notations, 
if $\phi_{\nu^-}$ is a positive definite function on $\Sigma_{\nu^-}$,
then all positive extensions $\phi_\nu$ of $\phi_{\nu^-}$ are in one-to-one
correspondence with contractions $\Ch_\nu:\Domain(\phi_{\nu^-}) 
\to\Range(\phi_{\nu^-})$.
\end{lemma}

\begin{proof} Note first that
\[
\big(\phi_{\nu^-}(s^{-1}t)\big)_{\{s,t\}\subset C_\nu, 
\{s,t\}\in E(\Ga_{\nu^-})}
\]
is a partially positive matrix whose corresponding pattern is the graph induced
by $\Ga_{\nu^-}$ on $C_\nu$. By applying Lemma~\ref{le:basic}, we obtain that
there exists a positive definite completion of this matrix, and, moreover,
that all completions are in one to one correspondence with
contractions $\Ch_\nu:\Domain(\phi_{\nu^-})\to \Range(\phi_{\nu^-})$.
Denote the completed
entry by $B(\Ch_\nu)$.

Now, any maximal clique in $\Gamma_\nu$ that is not a clique in 
$\Gamma_{\nu^-}$ is a translate of $C_\nu$, and the corresponding 
partially defined matrix detemined by $\phi_{\nu^-}$ is the same.
If we define then
\[
\phi_\nu(s)=
\begin{cases}
\phi_{\nu^-}(s)& \text{for } s\in\Sigma_{\nu^-},\\
B(\Ch_\nu)& \text{for } s=s_\nu,\\
B(\Ch_\nu)^* & \text{for } s=s_\nu^{-1}
\end{cases}
\] 
we obtain a positive definite function on $\Sigma_\nu$ that extends $\phi_{\nu^-}$.
It is easy to see that this correspondence is one-to-one.
\end{proof}

In particular, such extensions exist, and thus any positive definite
function on $\Sigma_{\nu^-}$ can be extended to a positive definite
function on $\Sigma_\nu$. Among them there exists a distinguished one, the 
\emph{central extension},
obtained by taking $\gamma_\nu=0$.

Now, combining Lemmas~\ref{le:obvious} and~\ref{le:step}, we obtain an 
extension theorem for positive definite functions.

\begin{theorem}\label{th:main}
If $\phi:\Sigma_\nu\to\LL(\HH)$ is a positive definite function, then
$\phi$ has a positive definite extension $\Phi:\FT\to\LL(\HH)$. 
There is a one-to-one correspondence between the set of all positive
extensions and the set of sequences of contractions $(\Ch_\mu)_{\nu\preceq\mu}$, where
$\Ch_\mu: \Domain(\Phi|\Sigma_{\mu^-}) \to \Range(\Phi|\Sigma_{\mu^-})$.
\end{theorem}

Note that in this correspondence $\Phi|\Sigma_{\mu}$ depends only on $\Ch_{\mu'}$
with $\nu\preceq \mu'\preceq \mu$, and thus the domain and
range of $\Ch_\mu$ are well defined by the previous $\Ch_{\mu'}$.
Then at step $\mu$ one can choose $\ga_\mu$ to be an arbitrary element
of the corresponding operator unit ball.
In particular, one obtains the \emph{central extension} by chosing at each
step $\ga_\mu=0$.

The most important consequence is an extension result for positive definite
functions on $S_n$. It should be noted that, contrary to $\Sigma_\nu$,
the sets $S_n$ do not depend on the particular lexicographic order~$\preceq$.

\begin{proposition}\label{pr:main}
Every positive definite function $\phi$ on $S_n$ has a positive definite extension
$\Phi$ on~$\FT$. If $\phi$ is radial (that is, $\phi(s)$ depends only on $|s|$), 
then one can choose also $\Phi$ radial.
\end{proposition}

\begin{proof}
The first part is an immediate consequence of Theorem~\ref{th:main},
obtained by taking $\Sigma_\nu=S_n$. Suppose $\Phi$ is positive definite; it is proved
in~\cite{FT}, 3.1 that,  if 
$\tilde\Phi(s)$ is the average of $\Phi$ on words of length $|s|$,
then $\tilde\Phi$ is also positive definite. The second part of the proposition
follows.
\end{proof}

Another consequence is the parametrization of all positive
definite functions. 

\begin{theorem}\label{th:main2}
There is a one-to-one correspondence between the set of all positive
definite functions $\Phi:\FT\to \LL(\HH)$ and the set of sequences of contractions $(\Ch_\mu)_{\mu\not=\{e\}}$, where
$\Ch_\mu: \Domain(\Phi|\Sigma_{\mu^-}) \to \Range(\Phi|\Sigma_{\mu^-})$.
\end{theorem}

The parametrization obtained is an analogue, for the free group,
 of the structure theory for classical Toeplitz matrices as presented,
 for instance, in~\cite{C}. We have chosen here to concentrate
 on the conceptual meaning, avoiding computational details. However, using~\eqref{eq:defect1},
 \eqref{eq:defect2} and~\eqref{eq:defect3},
  one can also write down precise formulas identifying the 
realization spaces appearing in~\eqref{eq:domrange} with 
  direct sums of defect spaces of some $\gamma_\mu$'s, with $\mu\preceq\mu^-$,
while repeated applications of~\eqref{eq:ugamma} yield unitaries connecting
  the different expressions for the same realization space.
  One can thus obtain explicit, but complicated formulas for 
  $\Domain(\Phi|\Sigma_{\mu^-})$ and $\Range(\Phi|\Sigma_{\mu^-})$,
  that would allow a recurrent definition
of the set of all families $ \Ch_\mu$ appearing in Theorem~\ref{th:main}, analogue to 
the \emph{choice sequences} of the classical case~\cite{FF}. 
  
  \smallskip

We may also rephrase Theorem~\ref{th:main} by using
Naimark's Theorem.

\begin{corollary}\label{co:rep}
Suppose $\phi:S_n\to\LL(\HH)$ is a positive definite function. Then
there exists a representation $\pi$ of\/ $\FT$ on a Hilbert space
$\KK$, and an operator $V:\HH\to\KK$, such that, for all
$s,t\in\FT$ such that $s^{-1}t$ has length $\le n$, we have
\[
\phi(s^{-1}t)=V^*\pi(s^{-1}t)V.
\]
\end{corollary}

There is an alternate way of looking at Theorem~\ref{th:main} that
will be useful in Section~\ref{se:noncom}.
Denote by $\MMM_n(\HH)$ the block operator matrices indexed by elements in
$S_n$; they have order $N(n)$ equal to the cardinality of $S_n$.
To any function $\phi$ on $S_{2n}$ corresponds
then a matrix $M(\phi)\in\MMM_n(\HH)$, defined by
$M(\phi)_{s,t}=\phi(s^{-1}t)$. We will call elements
$M=M_{s,t}\in \MMM_n$ whose entries depend only on $s^{-1}t$
\emph{Toeplitz}, and denote by $\TTT_n$
the space of Toeplitz matrices in $\MMM_n$.
It is spanned by the matrices
$\epsilon(\sigma)$, $\sigma\in S_n$, 
where $\epsilon(\sigma)_{s,t}=\delta_{s\sigma,t }$.

The above definition shows that
any $M_\phi$ is Toeplitz. The converse 
is also true, as shown by the next Lemma.

\begin{lemma}\label{le:toeplitz}
If $M\in \MMM_n$ is a Toeplitz matrix, then there exists 
a uniquely defined $\phi:S_{2n}\to\LL(\HH)$, such that $M=M_\phi$.
\end{lemma}

\begin{proof}
If $s\in S_{2n}$, we may write $s=t^{-1}r$, with $t,r\in S_n$. We define
then $\phi(s)=M_{t,r}$. The Toeplitz condition implies that this definition
does not depend on the particular decomposition of~$s$. In order to show
that $\phi$ is indeed positive definite, take $S\in \FT$, such that for
any $s,t\in S$, $s^{-1}t\in S_{2n}$. By adding elements 
to~$S$, if necessary, we may assume that the maximum of the
lengths $|s^{-1}t|$, for $s,t\in S$, is actually equal to $2n$.
Take then $s_0,t_0\in S$, such that $|s_0^{-1}t_0|=n$, and let 
$\hat s$ be the ``middle'' of the path (in the Cayley graph of $\FT$)
going from $s_0$ to $t_0$; that is, $\hat s$ is uniquely defined
by the conditions $|s_0^{-1}\hat s|=|\hat s^{-1}t_0|=n$. Using the tree
structure of the Cayley graph, it is then easy
to see that $|\hat s^{-1}s|\le n$ for any $s\in S$. So 
\[
A(\phi:S)=[\phi(t^{-1}r]_{t,r\in \hat s^{-1}S}
\]
is a submatrix of 
\[
M=[\phi(t^{-1}r]_{t,r\in S_n},
\]
and is thus positive definite. The uniqueness of $\phi$ is immediate.
\end{proof}

\begin{corollary}\label{co:rephrase}
If $T\in\TTT_n$, then
there exists a representation $\pi$ of $\FT$ and an operator $V$,
such that
\[
T_{s,t}=V^*\pi(s^{-1}t)V.
\]
\end{corollary}

\begin{remark}\label{re:2}
Theorems~\ref{th:main} and~\ref{th:main2} and their corollaries
extend to the case when $\mathbb{F}$ is a free
group with an infinite number of generators. The main
ideas are similar, but the details are more cumbersome, since
we have to use in several instances transfinite induction,
along the lines of~\cite{T}. 
\end{remark}

\section{The central extension}\label{se:central}

Suppose $\phi:S_n\to \LL(\HH)$ is positive definite. We have defined in 
Section~\ref{se:main},
for a fixed 
lexicographic order $\preceq$ on $\FT$,  the central extension of $\phi$ to
the whole of~$\FT$. We intend to prove in this section
that this central extension does not actually depend 
on the lexicographic order.

If $\Phi:\FT\to\LL(\HH)$ is positive definite, denote 
$\Phi_m=\Phi|S_m$. As in Lemma~\ref{le:obvious}, 
the map $\Phi\mapsto(\Phi|S_m)_{m\ge n}$ gives
a one-to-one correspondence between positive definite
extensions $\Phi$ of $\phi$ and sequences $(\Phi_m)_{m\ge n}$, with 
the property that $\Phi_n=\phi$, and, if $m\le m'$, then $\Phi_{m'}|S_m=\Phi_m$.
Then, in order
to show that the central extension does not depend on the order,
it is enough to show that $\Phi_{n+1}$ does not depend, since
then an induction argument finishes the proof.

If $\Phi$ is an extension of $\phi$, suppose
$(\KK_\Phi,\pi_\Phi,V_\Phi)$ is its minimal Naimark dilation. For $\Sigma\subset\FT$,
denote
\[
\KK_\Phi(\Sigma)=\bigvee_{\sigma\in\Sigma} \pi_\Phi(\sigma)V_\Phi\HH.
\]
Note, for further use, that for any $s\in\FT$, $\Sigma\subset\FT$
we have 
\[
\KK_\Phi(s\Sigma)=\pi_\Phi(s) \KK_\Phi(\Sigma).
\]

We will say that $\Phi$ is a \emph{maximal $n+1$ orthogonal} extension
of $\phi$ if the following property is satisfied: whenever $s,t\in\FT$,
$|s^{-1}t|=n+1$, and 
\begin{equation}\label{eq:sigma}
\Sigma=\{r\in\FT: \max(|r^{-1}s|, |r^{-1}t|)\le n\},
\end{equation}
we have
\begin{equation}\label{eq:maxort}
\big(\KK_\Phi(\Sigma\cup\{s\})\ominus \KK_\Phi(\Sigma)\big) \perp
\big(\KK_\Phi(\Sigma\cup\{t\})\ominus \KK_\Phi(\Sigma)\big).
\end{equation}
Note that it is enough to check this property for $s=e$.
Also, \eqref{eq:maxort} is equivalent to any of the equations
\begin{align}
\big(\KK_\Phi(\Sigma\cup\{s\})\ominus \KK_\Phi(\Sigma)\big)& \perp \big(\KK_\Phi(\Sigma\cup\{t\})\label{eq:maxort1};\\
\big(\KK_\Phi(\Sigma\cup\{t\})\ominus \KK_\Phi(\Sigma)\big)& \perp \big(\KK_\Phi(\Sigma\cup\{s\})\label{eq:maxort2}.
\end{align}

\begin{lemma}\label{le:maxort1}
If $\Phi,\Phi'$ are two maximal $n+1$ orthogonal extensions of $\phi$,
then $\Phi_{n+1}=\Phi'_{n+1}$.
\end{lemma}

\begin{proof}
Take $s\in\FT$, with $|s|=n+1$. We have to show that $\Phi(s)=\Phi'(s)$.

Denote by $(\KK,\pi,V)$, $(\KK',\pi',V')$ the minimal Naimark dilations of
$\Phi,\Phi'$ 
respectively. 
If
\[
\Sigma=\{r\in\FT : |r|\le n, |r^{-1}s|\le n\},
\] 
let $P_\Sigma, Q_e, Q_s$ be
the orthogonal projections (in $\KK$) onto $\KK(\Sigma)$, 
$\KK(\Sigma\cup\{e\})\ominus\KK(\Sigma)$, and 
$\KK(\Sigma\cup\{s\})\ominus\KK(\Sigma)$ respectively, and similarly  
$P'_\Sigma, Q'_e, Q'_s$. The maximal orthogonality assumption~\eqref{eq:maxort}
implies that $P_\Sigma, Q_e, Q_s$ are mutually orthogonal projections,
as well as $P'_\Sigma, Q'_e, Q'_s$.

By Corollary~\ref{co:chordal}, for
any $r,t\in  \Sigma\cup\{s\}$ we
have $|r^{-1}t|\le n$, and thus $\Phi(r^{-1}t)=
\Phi'(r^{-1}t)=\phi(r^{-1}t) $. Thus there is a unitary operator
$\Omega:\KK(\Sigma\cup\{s\})\to\KK'(\Sigma\cup\{s\})$, such that
$\Omega\pi(r)V=\pi'(r)V'$ for $r\in\Sigma\cup\{s\}$. 
A similar argument 
gives the existence of $\Xi:\KK(\Sigma\cup\{e\})\to\KK'(\Sigma\cup\{e\})$,
such that 
$\Xi\pi(r)V=\pi'(r)V'$ for $r\in\Sigma\cup\{e\}$.
moreover, $\Omega|\KK(\Sigma)=\Xi|\KK(\Sigma)$, and 
$\Omega P_\Sigma \Omega^*=P'_\Sigma,\ \Xi P_\Sigma \Xi^*=P'_\Sigma$.

Now, suppose $h\in\HH$. We have $V=(P_\Sigma+Q_e)V$ and
$\pi(s)V=(P_\Sigma+Q_s)\pi(s)V$, whence
\[
\begin{split}
\Phi(s)&=V^*\pi(s)Vh=V^*(P_\Sigma+Q_e)(P_\Sigma+Q_s)\pi(s)V
=V^*P_\Sigma\pi(s)V\\
&=V^*\Xi^*\Xi P_\Sigma \Xi^* \Omega P_\Sigma \Omega^* \Omega \pi(s)V
=V'{}^* P'_\Sigma \pi'(s) V'=\Phi'(s).
\end{split}
\]
The proof is finished.
\end{proof}

\begin{proposition}\label{pr:maxort2}
The central extension of $\phi$ corresponding to a given lexicographic order
is maximal $n+1$ orthogonal.
\end{proposition}

\begin{proof}
We have to prove~\eqref{eq:maxort} for all pairs $\{s,t\}$. 
We will use induction with respect to $\nu\in\hat\FT^*$ (note
that the lexicographic total order $\preceq$ is fixed): assuming that~\eqref{eq:maxort} is true whenever $\widehat{s^{-1}t}\prec \nu=\{r, r^{-1}\}$,
we will show that it is also true when $\widehat{s^{-1}t}=  \nu$. As
noted above, this is equivalent to show it for $s=e, t=r$.

Let us then fix, for the rest of this proof, $r\in\FT$, and assume 
$|r|=n+1$. Denote by $C$ the set of elements which are at distance
at most $n$ to both $e$ and $r$ and by $D$ the set of elements which are adjacent
to $e$ and $r$ in $\Gamma_{\nu^-}$; they are both cliques in 
$\Gamma_{\nu^-}$: $C$ by Corollary~\ref{co:chordal} (note that $\Gamma_n$
is a subgraph of $\Gamma_{\nu^-}$), and $D$ by  Lemma~\ref{le:not2cliques}
and Proposition~\ref{pr:2005}. We know by the definition of the central extension
with respect to the given lexicographic order that 
\begin{equation}\label{eq:ortD}
\big(\KK(D\cup\{e\})\ominus \KK(D)\big) \perp
\big(\KK(D\cup\{r\})\ominus \KK(D)\big),
\end{equation}
and we want to show that
\begin{equation}\label{eq:ortC}
\big(\KK(C\cup\{e\})\ominus \KK(C)\big) \perp
\big(\KK(C\cup\{r\})\ominus \KK(C)\big),
\end{equation}

Obviously $C\subset D$; and
for each $s\in D\setminus C$ we have either $|s|=n+1$ or
$|r^{-1}s|=n+1$. Let us denote
$ S_1=\{s\in D\cup\{e\} : |r^{-1}s|=n+1\}$, $ S_2=\{s\in D\cup\{r\} : |s|=n+1\} $.

We need another notation, also relative to the pair $\{e,r\}$.
For any $s\in\FT$, we define an element $\beta(s)$, as follows.
In the Cayley graph of $\FT$, suppose $P,P_1, P_2$ are the minimal paths connecting $e$ and $r$, 
$e$ and $s$, and $s$ and $r$ respectively. 
If $P\cap P_1=\{e\}$, we put $\beta(s)=e$, and if
$P \cap P_2=\{r\}$, then $\beta(s)=r$. (Note
that at most one of these equalities can be true,
since otherwise the union $P\cup P_1\cup P_2$
would be a cycle.) In the remaining case
the element on $P\cap P_1$ 
farthest from $e$ is also
the element on $P\cap P_2$ farthest from $r$; we will denote it 
by $\beta(s)$. One sees that $\beta:\FT\to P$, and $\beta(s)=s$ for
all $s\in P$.

Suppose $s\in D$; if $\beta(s)=e$, then $|r^{-1}s|>n+1$, which is 
a contradiction. Similarly one cannot have $\beta(s)=r$, and therefore
$\beta(s)\in P\setminus\{e,r\}$.

Let us now take $t\in C$, and $s\in S$. Suppose, for instance, that
$\beta(t)$ is closer to $r$ than $\beta(s)$. Since 
$|s^{-1}r|=|s^{-1}\beta(s)|+ |\beta(s)^{-1}r|$, 
$|r|=|\beta(s)|+|\beta(s)^{-1}r|$,
and $|s^{-1}r|\le n+1=|r|$, we have $|s^{-1}\beta(s)|\le |\beta(s)|$. Then
\[
|s^{-1}t|=|s^{-1}\beta(s)|+|\beta(s)^{-1}t|\le |\beta(s)|+|\beta(s)^{-1}t|=|t|=n.
\]
Therefore $t\in C$ implies that $|t^{-1}s|\le n$ for all $s\in D\cup \{e,r\}$.


%
%
%
%
%
%
%
%
%
%
%
%
%
%
%
%
%
%
%


Take now $s_1\in S_1$, $s_2\in S_2$.
Since
\[
|s_1^{-1}\beta(s_1)|+|\beta(s_1)^{-1}r|=|s_1^{-1}r|=n+1
=|r|=
|\beta(s_1)|+|\beta(s_1)^{-1}r|,
\]
we have $|s_1^{-1}\beta(s_1)|=|\beta(s_1)|$. Similarly,  $|s_2^{-1}\beta(s_2)|=|r^{-1}\beta(s_2)|$. Also, $\beta(s_1)$ is closer to $e$
than $\beta(r)$. Therefore
\[
\begin{split}
|s_1^{-1}s_2|&=|s_1^{-1}\beta(s_1)|+|\beta(s_1)^{-1} \beta(s_2)|+
|\beta(s_2)^{-1} r|\\
&= |\beta(s_1)|+|\beta(s_1)^{-1} \beta(s_2)|+|r^{-1}\beta(s_2)|
=|r|=n+1.
\end{split}
\]

Suppose now that there exists a vertex $t\notin C$, such that 
$|t^{-1}s_1|\le n$ and $|t^{-1}s_2|\le n$.
As shown above for elements of $D$, in this case also $\beta(t)\in P\{e,r\}$. 
If $\beta(s_1)$ is closer
to $e$ than $\beta(t)$, or if $\beta(s_1)=\beta(t)$ then
\[
\begin{split}
|t|&=|\beta(s_1)|+|\beta(s_1)^{-1}\beta(t)|+ |\beta(t)^{-1}t|\\
&= |s_1^{-1}\beta(s_1)|+|\beta(s_1)^{-1}\beta(t)|+|\beta(t)^{-1}t|
=|s_1^{-1}t|=n.
\end{split}
\]
In case $\beta(t)$ is closer to $e$ than $\beta(s_1)$, we obtain $|w|< n$.

Similarly, one proves that $|r^{-1}t|\le n$. Consequently, $t\in C$. It follows
then that 
\[
C=\{t\in \FT: d(t,s_1)\le n,\quad d(t,s_2)\le n\}.
\]
Now, if the pair $\{s_1, s_2\}$ is different from $\{e,r\}$, then, since
it is an edge of $\Gamma_{\nu^-}$, we must have $\widehat{s_1^{-1}s_2}=\mu$
for some $\mu\preceq\nu$.
By applying the induction hypothesis, it follows that
\begin{equation}\label{eq:ort1}
\big(\KK(C\cup\{s_1\})\ominus \KK(C)\big) \perp
\big(\KK(C\cup\{s_2\})\ominus \KK(C)\big).
\end{equation}
Therefore $(\KK( S_1\cup C)\ominus\KK(C))\perp (\KK( S_2\cup C)\ominus\KK(C))$,
and thus
\begin{equation}\label{eq:KKdec}
\KK(D)=\KK(C)\oplus 
\big(\KK( S_1\cup C)\ominus\KK(C)\big)\oplus
\big(\KK( S_2\cup C)\ominus\KK(C)\big)
\end{equation}
Since $e\in S_1$, \eqref{eq:ort1} implies in particular that $\KK(\{e\})\perp \KK( S_2\cup C)\ominus\KK(C)$. 
If $\xi\in \KK(\{e\})$, we can 
write the orthogonal sum
\[
\xi=\xi_\perp+\xi_C+\xi_1,
\]
with $\xi_\perp \perp \KK(D)$, $\xi_C\in\KK(C)$, and $\xi_1\in \KK( S_1\cup C)\ominus\KK(C)$.

Again by~\eqref{eq:ort1}, using the fact that $r\in S_2$, we have 
\[
\KK(\{r\})\perp \big(\KK( S_1\cup C)\ominus\KK(C)\big),
\]
and thus $\KK(\{r\})\perp \xi_1$. On the other hand, the centrality condition says
that $\xi_\perp\perp\KK(\{r\}) $. Therefore 
\[
\xi-\xi_C=\xi_\perp+\xi_1\perp\KK(\{r\}) ,
\]
which is exactly what have to prove.
\end{proof}

As a consequence of Lemma~\ref{le:maxort1} and Proposition~\ref{pr:maxort2}, we obtain
the main result of this section.

\begin{theorem}\label{th:maxort}
The central extension of $\phi:S_n\to\LL(\HH)$ does not depend on the lexicographic
order considered and is maximal orthogonal.
\end{theorem}

We may then speak about the \emph{central extension} of a positive definite
function defined on words of finite length, with no reference to a lexicographic
order on~$\FT$. Note that central liftings and central extensions have 
been extensively studied (see, for instance, \cite{BFF, FFG, GKW})

\section{Quasi-multiplicative functions}\label{se:quasim}

Haagerup has considered in~\cite{Ha} the functions
$s\mapsto e^{-t|s|}$, and has proved that they are positive definite
for $t>0$. In~\cite{MF} a larger class is defined: a function 
$u:\FT\to\bbC$ is called a \emph{Haagerup function} if $u(e)=1$,
$|u(s)|\le 1$, $u(s^{-1})=\overline{u(s)}$, and 
$u(st)=u(s)u(t)$ whenever $|st|=|s|+|t|$; it is proved therein
that any Haagerup function is positive definite. The result
is generalized in~\cite{Bo}, where Bo\.zejko introduces the analogue 
operator-valued functions (called \emph{quasi-multiplicative})
and proves that they are positive definite. Thus, we say that 
$\Phi:\FT\to\LL(\HH)$ is quasi-multiplicative if $\Phi(e)=I_\HH$,
$\|\Phi(s)\|\le 1$, $\Phi(s^{-1})=\Phi(s)^*$, and $\Phi(st)=\Phi(s)\Phi(t)$ whenever $|st|=|s|+|t|$.

One may say more about these functions in our context. We 
start with a preparatory lemma.

\begin{lemma}\label{le:qm}
Suppose $\Phi:\FT\to\LL(\HH)$ is quasi-multiplicative, 
and $(\KK,\pi,V)$ is
the minimal Naimark dilation of $\Phi$. If $|st|=|s|+|t|$,
then
\begin{equation}\label{eq:qm1}
\big(\KK(\{e,s\})\ominus \KK(\{s\})\big)\perp \KK(\{st\}).
\end{equation}
\end{lemma}

\begin{proof}
Take $\xi\in\KK(\{e\})$ and $\eta\in\KK(\{st\})$. By definition $\xi=Vh$
and $\eta=\pi(st)Vk$ for some $h,k\in\HH$, and we have
\begin{equation}\label{eq:xieta}
\<\xi, \eta\>=\<Vh,\pi(st)Vk\> = \<h,V^*\pi(st)Vk\>=\<h, \Phi(st)k\>.
\end{equation}

On the other hand, the orthogonal projection $P_s$ onto $\KK(\{s\})$ is
given by the formula $P_s=\pi(s)VV^*\pi(s)^*$. Therefore
\begin{align*}
\<P_s\xi, P_s\eta\>&= \<\pi(s)VV^*\pi(s)^*Vh , \pi(s)VV^*\pi(s)^*\pi(st)Vk \>\\
&=\<V^*\pi(s)^*Vh , V^*\pi(s)^*\pi(st)Vk \>
=\<V^*\pi(s)^*Vh , V^*\pi(t)Vk \>\\
&=
\< \Phi(s)^* h, \Phi(t)k\> 
=\< h,\Phi(s) \Phi(t)k\>=\<h, \Phi(st)k\>
\end{align*}
(we have applied quasi-multiplicity for the last equality).

Comparing the last relation with~\eqref{eq:xieta}, it follows that
$\<\xi, \eta\>=\<P_s\xi, P_s\eta\>$. Since
\[
\<\xi,\eta\>=\<P_{\KK(\{s\})}\xi,P_{\KK(\{s\})}\eta\> +
 \<\xi-P_{\KK(\{s\})}\xi,\eta\>,
\]
we obtain $\<\xi-P_{\KK(\{s\})}\xi,\eta\>=0$. But the space
on the left hand side of~\eqref{eq:qm1} is spanned by the vectors
$\xi-P_{\KK(\{s\})}\xi$, with $\xi\in\KK(\{e\})$;
thus the lemma is proved.
\end{proof}

One sees easily that a function $\phi:S_1\to\LL(\HH)$  with $\phi(e)=I$ is positive definite iff
$\phi(s^{-1})=\phi(s)^*$ and $|\phi(s)|\le 1$.

\begin{theorem}\label{th:bozejko}
A quasi-multiplicative function $\Phi:\FT\to\LL(\HH)$ is the
central extension of its restriction to $S_1$.
\end{theorem}

\begin{proof}
Denote by $\phi$ the restriction of $\Phi$ to $S_1$, and by $(\KK,\pi,V)$
the minimal Naimark dilation of $\Phi$.
According to Proposition~\ref{pr:maxort2}, we have to prove that $\Phi$
is maximal $n+1$ orthogonal for each $n\ge1$. 
 
Consider then $s,t\in\FT$, with
$|s^{-1}t|=n+1$, and suppose $\Sigma$ is defined by~\eqref{eq:sigma}. 
We will show that equality~\eqref{eq:maxort1} is true.

If $P$ is the minimal path connecting $s$ and $t$, suppose $x$ is the
element on $P$ between $s$ and $t$ adjacent to $s$. For any $r\in\Sigma\cup\{t\}$
we have $|s^{-1}r|=|s^{-1}x|+|x^{-1}r|$.
Then Lemma~\ref{le:qm} implies that
\[
\big(\KK(\{e,s^{-1}x\})\ominus \KK(\{s^{-1}x\})\big)\perp 
\KK(\{s^{-1}r\}).
\]
Since
$\big(\KK(\{s,x\})\ominus \KK(\{x\})\big)=\pi(s)
\big(\KK(\{e,s^{-1}x\})\ominus \KK(\{s^{-1}x\})\big)$ and
$\pi(s)\big( \KK(\{s^{-1}r\}) \big)= \KK(\{r\})$,
it follows that
\[
\big(\KK(\{s,x\})\ominus \KK(\{x\})\big) \perp \KK(\{r\}).
\]

Denote by $P_x$ the orthogonal projection onto $\KK(\{x\})$ and by $P_\Sigma$
the orthogonal projection onto $\KK(\Sigma)$.
The last
equality says that
for any $\xi\in\KK(\{s\})$ we have $\xi- P_x\xi \perp \KK(\{r\})$.
Therefore 
\begin{equation}\label{eq:perp}
\xi- P_x\xi \perp \KK(\Sigma\cup\{t\}). 
\end{equation}

In particular $\xi- P_x\xi \perp \KK(\Sigma)$; 
since obviously
$P_x\xi\in \KK(\Sigma)$, it follows that $P_x\xi=P_\Sigma\xi $.  But then~\eqref{eq:perp}
says that $\xi- P_\Sigma\xi \perp \KK(\Sigma\cup\{t\})$. Since the vectors of the form
$\xi- P_\Sigma\xi$, with $\xi\in \KK(\{s\})$, span $\KK(\{s\}\cup\Sigma)\ominus\KK(\Sigma)$,
relation~\eqref{eq:maxort1} follows; this ends the proof of the theorem.
\end{proof}

\section{Noncommutative factorization}\label{se:noncom}

We will apply in this section the above result to
noncommutative factorization problems, in the line of~\cite{MC}
and~\cite{H}. It is worth mentioning that the connection between
extension problems of positive definite functions and factorization
has been noted and used in the commutative case already in~\cite{Ru}.
The relation in the noncommutative case appears
in~\cite{MC} and~\cite{H}.
The technique used below is adapted from~\cite{MC}. A few
preliminaries are in order.

The elements of $\FT$ can also be considered as ``monomials'' in the indeterminates $X_1,
\dots, X_m$ and $X_1^{-1},\dots, X_m^{-1}$; the monomial $X(s)$
corresponding to $s\in\FT$ is obtained by replacing $a_i$ by $X_i$
and $a_i^{-1}$ by $X_i^{-1}$. It is then possible to
consider also polynomials in these indeterminates; that is, formal
finite sums
\begin{equation}\label{eq:def_p}
p(X)=\sum_{s} A_s X(s),
\end{equation}
where we  assume the coefficients $A_s$ to be operators on a fixed
Hilbert space $\CC$. We denote by $\deg(p)$ (the
\emph{degree} of $p$) the maximum length of the words appearing in
the sum in~\Ref{eq:def_p}. We introduce also an involution on
polynomials by defining, for $p$ as in~\Ref{eq:def_p},
$p(X)^*=\sum_{s} A_s^* X(s^{-1})$.

If $U_1,\dots U_m$ are unitary (not necessarily commuting) operators acting on a separable Hilbert space $\XX$,
then, for $s\in\FT$, $U(s)$ is the operator on $\XX$ obtained by replacing
$X_i$ with $U_i$, and $X_i^{-1}$ by $U_i^{-1}=U_i^*$.
Then $p(U)\in\LL(\CC\otimes\XX)$ is defined by $p(U)=\sum_{s} A_s \otimes U(s)$,
We say that such a polynomial is \emph{positive} if, for any
choice of unitary operators $U_1,\dots
U_m$, the operator $p(U)$ is positive.

\begin{theorem}\label{th:factorization}
If $p$ is a positive polynomial, then there exists an auxiliary
Hilbert space $\EE$ and operators $B_s:\CC\to\EE$, defined for
words $s$ of length $\le \deg(p)$ such that, if $q(X)=\sum_{s} B_s
X(s)$, then
\begin{equation}\label{eq:factorization}
p(X)=q(X)^* q(X)
\end{equation}
\end{theorem}

\begin{proof} As noted above, the proof follows the line of Theorem~0.1
in~\cite{MC}; the argument uses at a crucial point
Theorem~\ref{th:main}.

Denote $d=\deg(p)$, and consider the space of Toeplitz matrices $\TTT_d\subset\MMM_d$;
it is a finite dimensional operator system. We define a map
$\psi:\TTT_d\to\LL(\CC)$ by the formula
$\psi(\epsilon(\sigma))=A_\sigma$.

Suppose that an element $T\in\TTT_d \otimes \MM_k$ is a positive
matrix. By considering $T=[\tau_{s^{-1}t}]$ as an element of $\TTT_d(\MM_k)$, we can
apply  Corollary~\ref{co:rephrase} and obtain a representation
$\pi_\Phi:\FT\to\LL(\HHH_\Phi)$ and an operator $V_\Phi:\bbC^k\to \HHH_\Phi$, such
that
\[
\tau_{s^{-1}t}=V_\Phi^*\pi_\Phi(s^{-1}t)V_\Phi.
\]
Therefore
\begin{eqnarray*}
(\psi\otimes 1_k)(T)&=& (\psi\otimes 1_k) \left(\sum_s
\epsilon(s)\otimes \tau_s\right) = \sum_s  A_s\otimes \tau_s\\
& =& (1_\CC\otimes V_\Phi)^* \left(\sum_s A_h\otimes\pi_\Phi(s)
\right)(1_\CC\otimes V_\Phi).
\end{eqnarray*}

Since $p$ is positive, $\left(\sum_s A_h\otimes\pi_\Phi(s) \right)$
is a positive operator, and therefore the same is true about
$(\psi\otimes 1_k)(T)$. It follows then that $\psi$ is completely
positive.
Applying Arveson's Extension Theorem \cite{A, P}, we can extend $\psi$ to a
completely positive map $\tilde\psi:\MMM_n\to\LL(\CC)$.

Suppose $E_{s,t}\in\MMM_n$  has 1 in the $(s,t)$ position and 0
everywhere else. The block operator matrix $(E_{s,t})_{s,t}$
is positive; by 
Choi's Theorem~\cite[ch.~3]{P}
$(\tilde\psi(E_{s,t}))_{s,t}\in\LL(\bigoplus^N \CC)$ is positive,
where $N=N(d)$. Therefore, there exist
operators $B_s:\CC\to\bigoplus^N \CC$, such that
$B^*_sB_t=\tilde\psi(E_{s,t})$. Consequently
\begin{eqnarray*}
A_x&=& \psi(\epsilon(x))= \psi(\sum_{x=s^{-1}t} E_{s,t})\\
&=&\sum_{x=s^{-1}t} \tilde\psi (E_{s,t})=
\sum_{x=s^{-1}t}B^*_sB_t.
\end{eqnarray*}
If we define $q(X)=\sum_{s} B_s X(s)$, then the last equality is equivalent
to $p(X)=q(X)^* q(X)$.
\end{proof}

Note that the factorization of $p$ is usually not
unique. The theorem produces a factor $B$ which has degree at most
equal to $\deg p$. Also, in case $\CC$ is finite dimensional, the resulting
space $\EE$ can also be taken finite dimensional, with $\dim\EE=N(d)\times\dim\CC$.

One can rephrase the result of Theorem~\ref{th:factorization} as a
decomposition into sum of squares (making thus the connection
with~\cite{H}).

\begin{corollary}
If $p$ is a positive polynomial, then there exist a finite number
of polynomials $Q_1,\dots, Q_N$, with coefficients in $\LL(\CC)$,
such that
\begin{equation}\label{eq:sumofsquares}
p(X)=\sum_{j=1}^N Q_j(X)^* Q_j(X).
\end{equation}
\end{corollary}

\begin{proof} Since $B_s:\CC\to\bigoplus^N \CC$, we can write
$B_s=(B_s^{(1)}\ \dots\ B_s^{(N)})^t$. We consider then
$Q_j(X)=\sum_{s} B_s^{(j)} X(s)$, and obtain the required
decomposition.
\end{proof}

It is worth comparing the factorization obtained in
Theorem~\ref{th:factorization} with the results of~\cite{H}
and~\cite{MC}. In~\cite{H} the author considers real polynomials in the indeterminates
$X_1,\dots,X_m$ and $X_1^t,\dots,X_m^t$; $p$ is called positive
when any replacement of $X_i$ with real matrices $A_i$ and of $X_i^t$ with the transpose of $A_i$ leads to a positive semidefinite matrix. The main result is
an analogue of decomposition~\Ref{eq:sumofsquares} for such polynomials.

The analogue of Theorem~\ref{th:factorization} in~\cite{MC} deals with the case when the words appearing in the polynomial are of the form $s_+^{-1}t_+$, where $s_+, t_+$ are elements in the free semigroup with $m$ generators.
The positivity condition replaces the indeterminates with unitary matrices,
and the consequence is a corresponding decomposition.

\section*{Acknowledgements}

The authors are partially supported by NSF grant 12-21-11220-N65. The second author is partially supported by CERES grant 4-187 of the Romanian Government.


\begin{thebibliography}{xx}

\bibitem{ACC} \emph{Gr. Arsene, Z. Ceau\c sescu, T. Constantinescu}, Schur analysis of some completion problems. Linear Algebra Appl. {\bf 109}  (1988), 1--35.


\bibitem{A}
\emph{W.B. Arveson}, Subalgebras of $C^*$-algebras. Acta Math.
{\bf 123} (1969), 141--224.

\bibitem{B1} \emph{M. Bakonyi}, The extension of positive definite operator-valued
functions defined on a symmetric interval of an ordered group.
Proc. Amer. Math. Soc. {\bf 130} (2002), 1401--1406.

\bibitem{BC} \emph{M. Bakonyi, T. Constantinescu}, Inheritance principles
for chordal graphs. Linear Algebra Appl. {\bf 148} (1991), 125--143.

\bibitem{BN} \emph{M. Bakonyi, G. Naevdal}, 
The finite subsets of $\bbZ^2$ having the extension property.
J.~London Math. Soc. (2) {\bf 62} (2000), 904--916.


\bibitem{BFF} \emph{H. Bercovici, C. Foias, A. Frazho},  Central commutant liftings in the coupling approach. Libertas Math. {\bf 14}  (1994), 159--169. 

\bibitem{GR}  \emph{C. Berge}, The theory of graphs.  Dover Publications, Inc., Mineola, NY, 2001. 


\bibitem{Bo} \emph{M. Bo\.zejko},
Positive-definite kernels, length functions on groups and a
noncommutative von Neumann inequality. Studia Math. {\bf
95} (1989), 107--118.

\bibitem{BS} 
\emph{M. Bo\.zejko, R. Speicher}, 
Completely positive maps on Coxeter groups, deformed commutation relations, and operator spaces.
Math. Ann. {\bf 300} (1994), 97--120.

\bibitem{C} \emph{T. Constantinescu}, On the structure of positive Toeplitz
forms. Dilation Theory, Toeplitz Operators, and Other Topics,
127--149, Oper. Theory Adv. Appl., 11, Birkh\"auser, Basel, 1983.

\bibitem{C2} \emph{T. Constantinescu}, Schur analysis of positive block-matrices.
I.Schur Methods in Operator Theory and Signal Processing,
191--206, Oper. Theory Adv. Appl., 18, Birkh\"auser, Basel, 1986.

\bibitem{MF} \emph{L. De-Michele, A. Fig\`a-Talamanca},
Positive definite functions on free groups.
Amer. J. Math. {\bf 102} (1980), 503--509.

\bibitem{D} \emph{J. Dixmier}, Les $C^*$-alg\`{e}bres et leurs
repr\'{e}sentations.
Gauthier-Villars, Paris, 1969.

\bibitem{E} \emph{P. Eymard}, L'alg\`ebre de Fourier d'un groupe localement compact.
Bull. Soc. Math. France {\bf 92} (1964), 181--236.

\bibitem{FT} \emph{A. Fig\`a-Talamanca, M.A. Picardello}, Harmonic Analysis on Free
Groups. Marcel Dekker, New York, 1983.

\bibitem{FF} \emph{C. Foias, A.E. Frazho}, The Commutant Lifting Approach To 
Interpolation Problems. Operator Theory: Advances and Applications {\bf 44}, 
Birkh\"auser Verlag, Basel, 1990. 

\bibitem{FFG}  \emph{C. Foias, A.E. Frazho, I. Gohberg}, Central intertwining lifting, maximum entropy and their permanence. Integral Equations Operator Theory{\bf 18}  (1994), 166--201.

\bibitem{FFGK} \emph{C. Foias, A.E. Frazho, I. Gohberg, M.A. Kaashoek},
Metric constrained interpolation, commutant lifting and systems.
Operator Theory: Advances and Applications {\bf 100},
Birkh\"auser Verlag, Basel, 1998.

\bibitem{Ge} \emph{Ya. L. Geronimus}, Orthogonal Polynomials. Consultants Bureau, New York, 1961. 


\bibitem{God} \emph{R. Godement}, Les fonctions de type positif et la th\'{e}orie des
groupes. Trans. Amer. Math. Soc. {\bf 63} (1948), 1--84.

\bibitem{GKW}  \emph{I. Gohberg, M.A. Kaashoek, H.J. Woerdeman},  A maximum entropy principle in the general framework of the band method.
J. Funct. Anal.{\bf 95} (1991),  231--254.


\bibitem{GJSW} \emph{R. Grone, Ch.R.~Johnson, E.M.~de S\'{a}, H.~Wolkowicz},
Positive definite completions of partial Hermitian matrices.
Linear Algebra Appl. {\bf 58} (1984), 109--124.

\bibitem{Ha} \emph{U. Haagerup}, An example of a nonnuclear $C^{*} $-algebra,
which has the metric approximation property. Invent. Math.
{\bf 50} (1978/79), 279--293.


\bibitem{H} \emph{W.J. Helton}, ``Positive'' noncommutative polynomials are
sums of squares. Annals of Math. {\bf 156} (2002), 675--694.

\bibitem{HM} \emph{W.J. Helton, S.A. McCullough}, 
A Positivstellensatz for non-commutative polynomials. 
Trans. Amer. Math. Soc. {\bf 356} (2004), 3721--3737.

\bibitem{J} \emph{Ch.R. Johnson}, 
Matrix completion problems: a survey. 
Matrix theory and applications (Phoenix, AZ, 1989), 171--198,
Proc. Sympos. Appl. Math., {\bf 40},
Amer. Math. Soc., Providence, RI, 1990. 




\bibitem{K}
\emph{M.G. Krein}: Sur le probl\`eme de prolongement des functions hermitiennes
positives et continues. Dokl. Akad. Nauk. SSSR, {\bf 26} (1940), 17--22.


\bibitem{MC} \emph{S. McCullough}, Factorization of operator-valued polynomials
in several non-commuting variables. Linear Algebra Appl.
{\bf 326} (2001), 193--203.


\bibitem{P} \emph{V.I. Paulsen}, Completely Bounded Maps and Dilations. Longman,
New York, 1986.

\bibitem{Po0} \emph{Gelu Popescu}, Positive-definite functions on free semigroups.
Canad. J. Math. {\bf 48} (1996), 887--896.

\bibitem{Po1} \emph{Gelu Popescu},  Structure and entropy for Toeplitz kernels.
C. R. Acad. Sci. Paris Sér. I Math. {\bf 329} (1999),  129--134.

\bibitem{Po2} \emph{Gelu Popescu},  Structure and entropy for positive-definite Toeplitz kernels on free semigroups. J. Math. Anal. Appl. {\bf 254}  (2001),  191--218.


\bibitem{Ru} \emph{W. Rudin},  The extension problem for positive-definite functions.
Illinois J. Math. {\bf 7} (1963), 532--539.

\bibitem{Sa}
\emph{Z.~Sasv{\'a}ri},
Positive Definite and Definitizable Functions.
Akademie Verlag, Berlin, 1994.

\bibitem{T} \emph{D. Timotin}, Completions of matrices and the commutant lifting
theorem. J. Funct. Anal. {\bf 104} (1992), 291--298.




\end{thebibliography}
\end{document}